\documentclass[a4paper,UKenglish]{lipics}
 
\usepackage{microtype}

\title{To Teach Modal Logic: An Opinionated Survey \footnote{This work was partially supported by the project "Intentional Hybrid Logic" (MICINN FFI2013-47126-P); Department of Philosophy, Logic and Aesthetics; University of Salamanca (Spain). It has received also partial support from SeCyT (UNC-Argentina).}}
\titlerunning{To Teach Modal Logic} 

\author[]{Luis A. Urtubey}
\affil[]{Universidad Nacional de C\'{o}rdoba\\
  C\'{o}rdoba, Argentina\\
  \texttt{urtubey@ffyh.unc.edu.ar}}
  
\authorrunning{Luis\,A. Urtubey} 

\Copyright{Luis A. Urtubey}


\serieslogo{logo_ttl}
\volumeinfo
  {M. Antonia {Huertas}, Jo\~ao {Marcos}, Mar\'ia {Manzano}, Sophie {Pinchinat}, \\
  Fran\c{c}ois {Schwarzentruber}}
  {5}
  {4th International Conference on Tools for Teaching Logic}
  {1}
  {1}
  {243}
\EventShortName{TTL2015}

\begin{document}

\maketitle

\begin{abstract}
I aim to promote an alternative agenda for teaching modal logic chiefly inspired by the relationships between modal logic and philosophy. The guiding idea for this proposal is a reappraisal of the interest of modal logic in philosophy, which do not stem mainly from mathematical issues, but which is motivated by central problems of philosophy and language. I will point out some themes to start elaborating a guide for a more comprehensive approach to teach modal logic, and consider the contributions of dual-process theories in cognitive science, in order to explore a pedagogical framework for the proposed point of view.
\end{abstract}

\section{Introduction}
In his recent book \textit{Logic and how it gets that way}, Dale Jaquette \cite{jaquette} quotes a few lines of Nicholas Rescher's \textit{Autobiography}, where he talks about his conception of logic. Rescher's words read as follows:   
\begin{quote}
I never looked a logic as an accomplished fact --a field fully formed and merely in need of systematization and exposition, but saw it as an unfinished and malleable discipline to be reworked and structured for the sake of its applications. I always gravitated towards those logical issues and areas that were of potential utility for the elucidation and treatment of philosophical issues, and viewed logic not just as a self-contained discipline, but as a body of machinery for the accomplishment of extra-logical work. \cite{rescher} 
\end{quote}
I was struck by this quote of Rescher. Admittedly, it contains some insightful observations. Firstly, Rescher points out that logic can be seen as an unfinished discipline, which can be structured according to its applications. A second interesting point comes from the observation that philosopher's interest in logic is triggered by those issues which have at least potential utility for philosophical problems. Last but not least, logic is chiefly viewed  by Rescher as a tool for accomplishing extra-logical work. This last remark is relevant with respect to teaching logic, if we include it among those extra-logical activities. 

In this short paper, I will be concerned with these remarks from the perspective of modal logic in particular. It is clear that modal logic has been articulated by considering its applications, as it is witnessed by  the alternative interpretations of the modal operators of necessity and possibility. Thus, we are accustomed to talk nowadays about logical, metaphysical, deontic, and epistemic modalities, among many other options. Moreover, it is also easy to see that modal logic is central to the so-called Philosophical Logic, wherein those logical issues of philosophical concern are collected. Viewing modal logic as a tool for the accomplishment of extra-logical work is less frequent, but it is not totally unusual, as it is witnessed, for example, by certain applications in artificial intelligence. Our primary concern will be, in this respect the potential use of modal logic to improve people's reasoning abilities, by furnishing them with formal tools and normative standards. It has been said, however, that formal logic failed in its attempt to accomplish this objective. Consequently, modal logic  could not reach this honour either. Contrariwise, I will argue that assuming and elaborating certain perspective with respect to cognition and reasoning, modal logic can still be fruitfully applied to improve people's reasoning abilities. Notwithstanding, it should not be misinterpreted. I will not offer here any specific details about the concrete organization of a modal logic course following this approach. I am only giving an overview, which can be exploited and implemented in different ways.

\section{Modal Logic: Philosophical Issues}

I am going to take a look now at some themes which populated the relationships between modal logic and philosophy during the twentieth century. I think that these examples can serve to articulate a bunch of topics for teaching modal logic. Admittedly, this is not an exhaustive enumeration, but it will suffice to give an idea of the approach we have in mind. I will articulate these issues in a pedagogical approach towards the end of the paper. 
A very useful summary to start with  is given by the comprehensive article by Linsdstr\"{o}m and Segerber in \textit{The Handbook of Modal Logic} \cite{linsdst-segerb}, where they point out some themes concerning the relationship between modal logic and philosophy that are also relevant from our point of view. They enumerates the followings:
\begin{enumerate}
\item
The back-and-forth between philosophy and modal logic: good examples are Carnap and Prior among many others, who tried to use  modal logic to throw light on old philosophical questions.
\item
The interpretation problem. The problem of providing a certain modal logic with an intuitive interpretation, which  should not be conflated with the problem of providing a formal system with a model-theoretic semantics.
\item
The fact that one may compare this situation with that in probability theory, where definitions of concepts like ‘outcome space’ and ‘random variable’ are orthogonal to questions about interpretations of the concept of probability.
\item
The value of formalization. Modal logic sets standards of precision, which are a challenge to —-and sometimes a model for-— philosophy. Classical philosophical questions can be sharpened and seen from a new perspective when formulated in the framework of modal logic. On the other hand, representing old questions in a formal garb has its dangers, such as simplification and distortion.
\item
Why modal logic rather than classical (first or higher order) logic? The idioms of modal logic seem better to correspond to human ways of thinking than ordinary extensional logic. One can think, for example about counterfactual reasoning and its central place in human thinking. 
\end{enumerate}
Let's go a little bit more along this path.

\subsection{The Readings of Modality}

As it was said above, Carnap and Prior constitute two important philosophical landmarks in the development of modal logic. Thence, it is worthwhile to isolate the philosophical approach to modalities of these authors and to introduce thereby other themes which are related with their work. Lindstr\"{o}m and Segerberg observe that Carnap's project was not only to develop a semantics (in the sense of Tarski) for intensional languages, but also to use metalinguistic notions from formal semantics to throw light on the modal ones. In "Modalities and quantification" from 1946, he writes:
\begin{quote}
It seems to me ... that it is not possible to construct a satisfactory system before the meaning of the modalities are sufficiently clarified. I further believe that this clarification can best be achieved by correlating each of the modal concepts with a corresponding semantical concept (for example, necessity with L-truth).
\end{quote}
Concerning Prior, years ago, in a remarkably thoughtful article \cite{menzel}, Christopher Menzel analyzed Prior's approach to modal logic and his deep philosophical concerns. Menzel points out that, for Arthur Prior, "the construction of a logic was a supremely philosophical task". Being a clever logician, Prior joined the virtue of appreciating "a finely crafted formal system for its own sake", apart of its "meaning" or philosophical significance, as well as the power of a genuine logic "to lays bare the nature of those concepts that it purports to be a logic of". This last task can be accomplished only by way of deep reflection and insightful philosophical analysis.
Menzel remarks that Prior's  search for the “true” modal logic is where this sort of reflection and analysis is more evident\footnote{Particularly Menzel considers system Q from Prior's "Tense Logic for Non-Permanent Existents"\cite{prior}}.  Menzel's paper in its turn is a nice guide "to trace the course of Prior's own struggles with the concepts and phenomena of modality, and the reasoning that led him to his own rather peculiar modal logic Q". It is also an advisable starting point to appreciate Prior's intuitions and the arguments that rest upon them.

Let us also consider Prior's work under a new light. In a recent article about Arthur Prior's work on hybrid logic \cite{blackburn}, Patrick Blackburn has considered extensively those aspects of Prior's philosophy which motivated his interest on hybrid logic. This kind of logic departs from traditional logical categories, allowing that formulas can be used as terms. Blackburn's article is also an excellent starting point to stimulate the reflection on philosophical questions, which have motivated the development of modal languages. Blackburn points out that Prior's development of hybrid logic was fuelled by “his most fundamental convictions concerning time and tense, namely that the modal perspective, and in particular the internal view offered by modal logic, was the way to capture genuinely temporal logic”. In the early years of the last Century, the Hegelian philosopher John McTaggard posed a tricky argument concerning the unreality of time. McTaggard's argument is based on the postulation of two alternative series, which respectively constitute two ways of conceptualizing time: the A series and B series conceptions. In his work \textit{Past, Present and Future}, Prior discusses McTaggard argument. The A series considers the flow of time from past, through present, to future. It is an internal and situated conception of time, which reflects the experience of temporality of human beings.

Under the B series, time is a collection of instants ordered by a relation of temporal precedence, i.e. positions are ordered from earlier-than to later-than relations. Particularly, Prior favours A series talk over B series talk, because he thought that the external perspective underlying B series did not reflect the way in which human beings experience time. This motivate, in its turn the development of hybrid logic, since Prior believed that hybrid logic was the A series logic underpinning everyday tensed logic. 
There are too many connections between Prior's temporal logic and philosophical problems raised by Blackburn's article. I cannot survey all of them here. Nonetheless, the following quotations from \cite{blackburn} will be useful to illustrate an important point concerning the changes experimented by the relationship between logic and philosophy through the years:
\begin{quotation}
It is important to realize that Prior's willingness, so clearly in evidence here, to base the core of a philosophical position at least on the differing properties of formal logical systems, contrasts sharply with the attitude of many, perhaps most, contemporary modal logicians. Nowadays logicians are rather more guarded about the direct relevance of the properties of particular formal systems to philosophy. (…)  Much of this contemporary attitude stems from the model-theoretic perspective which underlies current practice.   
\end{quotation}
Nonetheless, Patrick Blackburn also notices that
\begin{quote}
(...) this does not mean that the model-theoretic conception rules out the possibility of drawing philosophical conclusions based on logical analysis.
\end{quote}
On the contrary: 
\begin{quote}
(...) The most direct way that model-theoretically oriented logic helps philosophy is by demonstrating that there are more ontological possibilities than is apparent at first sight, or that certain possibilities can't be made to fit together coherently. The model-theoretic perspective helps induce philosophical modesty.
\end{quote}

This also holds for Prior's work on tense and temporal logic.

Admittedly, unlike the case of mathematical modal logic, for teaching modal logic in philosophy or broadly speaking in the humanities, it seems to be crucial to count with a substantial reading of the box or the diamond. Following ideas of an unpublished paper by John P. Burgess \cite{burgess}, the point to be observed is that we have to admit throughout the whole history two alternative readings of the box as "it is analytic that...", which seems the official explanation of I. Lewis, and reading the box as "it could not have failed to be the case that..." (a reading which doubtless influenced the "intuitions" of some). It must also  be admitted that these were two different and incompatible readings. We'll go back on this below. 

\subsection{Quine's Interpretational Challenge}

For the second item above mentioned, it is inevitable to think of  Quine's criticism of quantified modal logic\footnote{The most known references at this respect are \cite{quine1}, \cite{quine2}, \cite{quine3}, \cite{quine4}}. According to the article of Lindstr\"{o}m and Segerberg, this criticism comes in different strands. First, there is the simple observation that classical quantification theory with identity cannot be applied to a language in which substitutivity of identicals for singular terms fails. The so-called Morning Star Paradox shows that we must abandon either universal specification (and existential generalisation) or indiscernibility of identicals. This premise goes to the seemingly  uncontroversial conclusion that classical quantification theory (with identity and individual constants) cannot be mixed up with non-extensional operators (those for which substitutivity of identicals for singular terms fail) without being modified in some way. This weak claim already gives rise to the challenge of extending quantification theory in a consistent way to languages with non-extensional operators.\footnote{A very comprehensive treatment is given by James Garson in \cite{garson2}} 
Perspicuously Lindstr\"{o}m and Segerberg differentiate a  much stronger claim in addition to this former claim, which sometimes can be found in Quine's early works, relative to the fact that objectual quantification into non-extensional (“opaque”) constructions simply does not make sense. The cause is that occurrences of variables inside of opaque constructions do not serve simply to refer to their objects, and cannot therefore be bound by quantifiers outside of the opaque construction. Consequently, quantifying into opaque contexts  would be like trying to quantify into quotations. As Lindstr\"{o}m and Segerberg observe: "This claim is hardly credible in the face of the multitude of quantified intensional logics that have been developed since then, and it can be taken to be refuted by the work of Kaplan and Fine among others.
I think that Lindstr\"{o}m and Segerberg's distinction about Quine's very different thesis is a fact that must be appreciated and incorporated when teaching modal logic in order to exploit this frequent confusions.

\subsection{The Nature of Modality}

It will be worthwhile to return here to some reflections about modal logic recently introduced by John P. Burgess that we just referred to \cite{burgess}. According to Burgess' approach, it is convenient to look backward to the origins of modal logic. Originally logic  used to "name an enterprise aiming to prescribe the norms of deductive argumentation in philosophy and elsewhere". However, "the pursuit of what were the original aims of logic is now for many inseparable from philosophy generally, and philosophy of language especially.\cite{burgess} Thus, logic progressively has gone far away from its original motivations. In fact, logic is no more restricted to the study of correct or valid reasoning. Formal systems are studied in its own right. Burgess' observation also enables to see that the original logic's aim is now shared by many disciplines other than logic. In consequence, in order to fulfil this original aim one must cross some disciplinary borders, since it is no longer in the domains of logic alone. Burgess also considers the case of modal logic in particular and thinks that we must go back to the work of Aristotle to appreciate the initial sense of modality that was getting lost as time went by. 
Moreover, a neglected question at this point is what is meant by "modal" here. And if there is a sense which is particular relevant to philosophy. What is characteristic of modality is the expression of irreality of one kind or another, as it is expressed by means of grammar. It is widely accepted the distinction between deontic, epistemic and dynamic modalities. In the beginnings of modal logic, H. von Wright also made an almost similar distinction. In fact, all of them were of interest in philosophy, but specially the third, which more or less corresponds to what is called "metaphysical" modality. This is the conclusion of Burgess, which is worth to consider in the setting of teaching modal logic in philosophy. The notions associated with this type of modality are the necessary, in the sense, of what is and would have been (could not failed to be) no matter what, and the possible in the sense of what is or isn't but potentially might (or could) have been. 
The idea of a true modal logic endorsed by Burgess and other philosophers must be interpreted as an effort to put modal logic again in the right way. In the words of Burgess \cite{burgess}:
\begin{quote}
Very little of modal logic is devoted to the attempt to determine what are the valid forms of argument involving necessity and possibility in any sense, and vanishingly little is devoted to the original aim of modal logic, that of dealing with the logic of modality in our default sense.
\end{quote}
The starting point in this regard is a distinction, which admits two contrasting opinions about modality that have coexisted through the whole history of modal logic.
Together with this initial idea of necessity, there would be an idea of necessity as analyticity, which is also unconsciously influenced by thoughts about necessity in the "metaphysical" sense. According to Burgess, the general pre-Kripkean failure to recognize both notions explicitly and to distinguish between them naturally made it impossible to appreciate that different logics may be appropriate to the different notions. As it happens with the precedent observations in this paper, I believe that this thesis about the confusion between two alternatives ideas of necessity can be usefully incorporated in a strategy for teaching modal logic.
Among other confusions that one can take into account regarding the application of modal logic, Burgess also considers "the modal distinction between the contingent and the necessary". Let us consider the following examples also from \cite{burgess}.

\begin{enumerate}
\item 
There has been no female U. S. president.
\item 
The number 29 has no nontrivial divisors.
\end{enumerate}

Admittedly, though there has been no female U. S. president, there could have been, whereas the number 29 not only has no nontrivial divisors but couldn't have had any. Thence, it is allowed to make counterfactual speculations regarding the first proposition, but not with respect to the second. And the most important question is how we know about the difference between these propositions. This is a central problem for modal epistemology concerning the bridge between 'is' and 'could'. It is also an interesting problem to reflect upon when you teach modal logic. A nice place to start is the classic statement of the problem given by Immanuel Kant in the introduction to the second or B edition of his \textit{Critique of Pure Reason} that Burgess points out:
\begin{quote}
Experience teaches us that a thing is so, but not that cannot be otherwise.
\end{quote}
That is to say that experience can teach us that a necessary truth is true, but it cannot teach us that it is necessary. Translated into the lexicon of possible world semantics it means that we can access through our senses to the actual world -the one we live in- whereas when we speak about what could or couldn't have been, we are making an assertion about how thinks are in other worlds. 
\subsection{A Modal Puzzle}
\label{sect:a modal puzzle}
Let us close this section by considering a well-known tricky puzzle about modality which can swell the list of teaching devices. The question concerning the right modal logic for "metaphysical modality" has already prompted the interest of philosophers and logicians. The most favoured systems have been S5 and S4, though there have been alternatives like McKinsey S4.1. Nathan Salmon and others have developed an argument for the system T, or rather against S4. As you know the characteristic law distinguishing S4 and T claims that if a proposition P is necessary, it is necessarily necessary (or if P is possibly possible), then it is possible. Salmon's argument has became to be known as the "modal paradox". Let us now describe the example which centres in the S4 axiom following the exposition of Burgess \cite{burgess}.
\begin{quote}
The same principle [S4] would seem to apply to, say, a ship made of a thousand planks. There is an intuition X to the effect that the same ship could have been made of 999 of the same planks, arranged the same way, plus a replacement for plank 473; but at the same time there is an intuition Y to the effect that a ship made of a thousand different planks would have been a different ship. The puzzle is that these two intuitions, X and Y, seem to conflict. 
\end{quote}

In \cite{burgess}, the conflict referred above is described in Leibnizian-Kripkean possible worlds language as follows.

\begin{quote}
Here in the actual world A we have a ship, let us name it the good ship Theseus, made of 1000 planks. Our first intuition X is that the same ship could have been made of 999 of these planks plus a replacement for plank 473. In possible-worlds terms that means there is another world B where the same good ship Theseus exists with all but one plank the same as in our world A, and only plank 473 different. (...) [A]nd by the same sort of intuition, there must be another world C where the same good ship Theseus exists with all (...) but with plank 692 different. That means for us back in world A there is another world C where the good ship Theseus exists with all but two planks the same as in our world A, but with planks 473 and 692 different (...). The same sort of considerations can then be used to argue that one could have three planks different, or four, or five, or all 1000. But that is contrary to our other intuition Y. 
\end{quote}

Salmon's solution rejects S4 axiom to avoid treating identity as vague. It is tantamount to recognize that though world B is a possibility for us in world A, and world C is a possibility for those in world B, still C is not a possibility for us in world A, but only a "possible possibility". That means that transitiveness must be banned in this context. There are reasons, however, as Burgess claims, to doubt of this analysis. Since the puzzle seems to be parallel to another famous temporal one already discussed by Plutarch and Thomas Hobbes, which involves the ship of Theseus. And in this last case, transitivity cannot be excluded.   
It is admissible to conjecture that this problem concerning the necessity of origins will have great significance and may ultimately be relevant to discuss the status of the law that possibly possible implies possible, when we have to teach about normal systems of modal logic. 

\section{Dual-Process Theories}
\label{sect:dual--process theories}
Dual-process theory in cognitive psychology is closely related to the study of heuristics and biases in reasoning, as it was taken up for example in the research of Kahneman and Tversky \cite{kahneman}. From this point of view, heuristics have intrinsic limits. And it is claimed that heuristics inserted in the human mind by the process of evolution can produce errors, although it is also given to human beings to fix these failures. By means of different tasks, Stanovich together with other researchers have shown the conflict between heuristics and  the so-called "analytic system". Precisely, the goal of these tasks is to bring about "a heuristically triggered response against a normative response generated by the analytic system" \cite{stanovich-toplak-West-2008}. The idea of an "analytic system" is crucial: it is one of the most important contributions of Kahneman and Tversky that cognitive biases are not made by chance, since the human mind utilizes "innate rules of thumb", which are autonomous and automatics, and which cannot be totally banned by us when processing information. Moreover, Kahneman and Frederick also have remarked that "[t]he persistence of such systematic errors in the intuitions of experts implied that their intuitive judgements may be governed by fundamentally different processes than the slower, more deliberate computations they had been trained to execute" \cite{kahneman-frederick-2005}.
It is then possible to design certain tasks or contexts where these natural heuristics bring out responses which are not adequate, before being corrected by means of a second system that is slower, more analytic, and more thoughtful. Thence, it is one of the central intuitions of dual-process theories that the human mind is composed by two different systems, the so-called system1 (S1) and system2 (S2). These models are heavily discussed in the recent literature of cognitive science though. The collection of processes of S1 (the modules) are autonomous and cognitively closed; we are not conscious of their working and we cannot influence on them directly. These processes are experimented in our intuitions and our spontaneous reactions. The S2 processes, on the contrary, are voluntary and they are able to inhibit and to fix those automatic responses of S1. Most of the time, these two systems get on well with each other, but there are situations where the response elicited by S1 conflicts with that issued by S2. Nevertheless, we can learn the rules and we can learn to apply them rightly by using S2 processes. According to the advocates of dual-process theories, for accomplishing a task it is necessary to identify first the bias of S1 and then to fix it with the appropriate tools of S2. Moreover, we have a tendency to utilize only S1, which in most cases gives us right answers, but the processes of this system mislead us in other cases. Fortunately, we may learn how to utilize S2  to detect these situations, where we are cheated, in order to look for the normatively adequate answer.
Stanovich \cite{Stanovich-2009} has stressed the importance of acquiring the cognitive tools which make it possible to attain the resolution of those problems which S1 has not answered adequately. These are the kind of tools which are usually provided by logic courses. However, as it has been emphasized, this is not enough to overcome this obstacle: one has to have "the relevant logical rules and also [needs] to recognize the applicability of these rules in particular situations" \cite{kahneman-frederick-2002}. Thence it must be also detected that the automatic answer provided by S1 has to be neutralized. Consequently biases' correction is more complex and at least involves two stages. This is the most important thing to be learned from these studies concerning their pedagogical application. Particularly, the complexity of this dual process serves also to explain the shouting failure of standard logic courses with respect to improve the reasoning abilities of many students.
From this point of view, Guillaume Beaulac and Serge Robert \cite{beaulac-robert} have pointed out that we can rescue the normative role of logic, by considering it in the light of a pedagogical reading of dual-process theories. The idea can be applied also to modal logic. To that end efforts must be made to elaborate cases which prompted the action of S1-systems that will have as normative counterparts the current systems of modal logic. In this case, S1-systems will serve to show up several drawbacks of modal reasoning and other puzzling issues which stem from reasoning tasks performed in this context. I think that the philosophical problems that have been shortly addressed in this paper, can --at least partially-- fulfil the requirements to confront the intuitions given by S1-systems in this case and contribute to improve teaching modal logic under these premises. 
  
\section{Conclusion}
\label{sect:conclusion}
In this paper, I aimed at integrating from a pedagogical perspective, the formal treatment of modalities with those philosophical problems which appeared in the origin of modal logic and have also contributed to its development. I have considered modal logic as a tool for treating and elaborating many problems in philosophy and other connected areas, relegating the mathematical issues. Concerning teaching, the dose of mathematical modal logic to be introduced must be estimated from this point of view and it will no more be the only motive of concern. On the contrary, it turns out to be instrumental. Specifically I have not considered  the  mathematics of modal logic, which I presupposed in a certain way. The central tenet may be summarized by saying that things must be turned upside down. Instead of focusing on the formal and mathematical aspects of modal logic, it is convenient to return to the origins and search firstly for the philosophical motivation behind the primeval ideas.

More controversially, the strategy for teaching modal logic that I am proposing here would not be that new, but I confess that I have never seen a full-fledged development of this approach. Certainly, one can find some exemplifications or applications of modal logic to philosophical thesis scattered in different logic texts.  I have surveyed in this paper some of the most relevant and more accessible ones in order to contribute to the construction of an alternative perspective for teaching modal logic in philosophy.  It may sound perhaps too much oriented to philosophy ignoring other uses of modal logic. That might be right, but I think that for teaching modal logic in philosophy it must be taken into account the central motivations of those people who have decided to study philosophy. And the challenge for all of us involved both in modal logic and philosophy is to articulate courses, which pay attention to these particular motivations. 

In this regard, our position lies close to some recent trends in the development of Mathematical Philosophy, which offers a new approach to philosophical problems chiefly supported  by the contribution of formal methods. The most innovative feature of this proposal is the balanced combination of formal methods with deep philosophical questions, which  trigger a productive reflection. This relationship suggests  also a link between this approach and the application of dual-process theories to the extra-logical activity of teaching and learning the application of formal tools.

\subparagraph*{Acknowledgements}
The author gratefully thank two anonymous referees of this publication for their useful comments. 



\bibliography{dummybib}



\end{document}